\documentclass[12pt,reqno]{amsart}
\usepackage{amsmath, amsfonts,epsfig, amssymb,amsbsy,eucal,mathrsfs,float, amsthm,mathtools,amsrefs}

\usepackage[normalem]{ulem}

\usepackage{graphicx,xcolor}
\usepackage{dynkin-diagrams}
\usepackage{fullpage}
\usepackage{textpos}

\usepackage{tikz}

\usepackage[all]{xy}

\usepackage[enableskew]{youngtab}

\usepackage{tikz-cd}
\usepackage{tikz}

\usepackage{hyperref}
\hypersetup{
colorlinks = true,
urlcolor = blue,
citecolor = .
}

\usepackage{todonotes}

\numberwithin{equation}{section}
\theoremstyle{plain}

\newtheorem{theorem}{Theorem}[section]
\newtheorem{lemma}[theorem]{Lemma}

\newtheorem{conjecture}[theorem]{Conjecture}

\theoremstyle{definition}

\newtheorem{definition}[theorem]{Definition}
\newtheorem{example}[theorem]{Example}

\theoremstyle{remark}

\newtheorem{remark}[theorem]{Remark}

\DeclareMathOperator{\G}{Gr}
\DeclareMathOperator{\Q}{Q}
\DeclareMathOperator{\IG}{IG}
\DeclareMathOperator{\OG}{OG}

\DeclareMathOperator{\Sp}{Sp}
\DeclareMathOperator{\QH}{QH}

\DeclareMathOperator{\BQH}{BQH}
\DeclareMathOperator{\Ext}{Ext}
\DeclareMathOperator{\Hom}{Hom}

\DeclareMathOperator{\Aut}{Aut}

\DeclareMathOperator{\Spec}{Spec}

\DeclareMathOperator{\SO}{SO}

\newcommand{\bA}{{\mathbb A}}
\newcommand{\bC}{{\mathbb C}}

\newcommand{\bS}{{\mathbb S}}

\newcommand{\bP}{{\mathbb P}}

\newcommand{\rA}{\mathrm{A}}
\newcommand{\rB}{\mathrm{B}}
\newcommand{\rC}{\mathrm{C}}
\newcommand{\rD}{\mathrm{D}}
\newcommand{\rE}{\mathrm{E}}
\newcommand{\rF}{\mathrm{F}}
\newcommand{\rG}{\mathrm{G}}

\newcommand{\rP}{\mathrm{P}}
\newcommand{\rT}{\mathrm{T}}
\newcommand{\rTs}{\mathrm{T}_{\mathrm{short}}}

\newcommand{\cO}{\mathcal{O}}
\newcommand{\cR}{\mathcal{R}}

\newcommand{\cU}{\mathcal{U}}

\newcommand{\Db}{{\mathbf D^{\mathrm{b}}}}

\renewcommand{\emph}{\textsf}

\DeclareMathOperator{\QS}{QS}
\newcommand{\QSo}{\QS^\circ}
\newcommand{\QSx}{\QS^\times}
\newcommand{\QHcan}{\QH_{\mathrm{can}}}

\begin{document}

\setcounter{section}{1}


\title{Residual categories of Grassmannians}

\author{Maxim N. Smirnov}
\address{
\parbox{0.95\textwidth}{
Universit\"at Augsburg,
Institut f\"ur Mathematik,
Universit\"atsstr.~14,
86159 Augsburg,
Germany
\smallskip
}}
\email{maxim.smirnov@math.uni-augsburg.de}
\email{maxim.n.smirnov@gmail.com}

\thanks{M.S. was partially supported by the Deutsche Forschungsgemeinschaft (DFG, German Research Foundation) --- Projektnummer 448537907.}

\maketitle

\begin{abstract}
This short note is an extended abstract for my talk at the
Nottingham Online Algebraic Geometry Seminar on October 1, 2020.
It is based on the joint works with Alexander Kuznetsov \cites{KuSm20, KuSm21}.
\end{abstract}


\bigskip

We work over a fixed algebraically closed field $\Bbbk$ of characteristic zero.

\subsection{Exceptional collections}

Let $X$ be a smooth projective variety over the field $\Bbbk$ and let $\Db(X)$ be
the bounded derived category of coherent sheaves on $X$. An object $E$ of $\Db(X)$
is called \emph{exceptional}, if we have
\begin{equation*}
  \Hom(E,E) = \Bbbk id_E \quad  \text{and} \quad \Ext^i(E,E) = 0 \quad  \forall i \neq 0.
\end{equation*}
A sequence of exceptional objects $E_1, \dots, E_n$ is called an \emph{exceptional collection},
if we have
\begin{equation*}
\Ext^k(E_i,E_j) = 0 \quad \text{for $i > j$ and   } \forall k.
\end{equation*}
We denote by $\langle E_1, \dots, E_n \rangle$
the smallest full triangulated subcategory of $\Db(X)$ containing the objects
$E_1, \dots, E_n$. If we have
\begin{equation*}
  \Db(X) = \langle E_1, \dots, E_n \rangle,
\end{equation*}
then the collection $E_1, \dots, E_n$ is called \emph{full}.

\begin{example}
  \
  \begin{enumerate}
    \item On the projective space $\bP^n$ we have the famous Beilinson collection \cite{Be}
    \begin{equation}\label{eq:collection-beilinson}
      \Db(\bP^n) = \langle \cO_{\bP^n}, \cO_{\bP^n}(1), \dots, \cO_{\bP^n}(n) \rangle.
    \end{equation}

    \item For Grassmannians $\G(k,n)$ and quadrics $\Q_n$ full exceptional collections
    were constructed by Kapranov \cite{Ka}. In the case of $\G(2,4)$ Kapranov's collection
    takes the form
    \begin{equation}\label{eq:collection-kapranov}
      \Db(\G(2,4)) = \langle \cO, \cU^\vee, S^2 \cU^\vee, \cO(1), \cU^\vee(1), \cO(2) \rangle,
    \end{equation}
    where $\cU$ is the tautological subbundle on $\G(2,4)$.
  \end{enumerate}
\end{example}

\subsection{Lefschetz exceptional collections}

In this talk we are interested in a particular class of exceptional collections,
called \emph{Lefschetz collections}
introduced by Alexander Kuznetsov in~\cite{K07} in the context of
homological projective duality. The main goal of the talk is to explain that
Lefschetz collections also have close connections with quantum cohomology and
mirror symmetry.

\begin{definition}
  Let $X$ be a smooth projective variety over $\Bbbk$ and let $\cO(1)$ be a line
  bundle on $X$. For an object $F \in \Db(X)$ we denote $F(1) \coloneqq F \otimes \cO(1)$.
  \begin{enumerate}
    \item[(i)] A \emph{Lefschetz collection} with respect to $\cO(1)$ is an exceptional collection,
    which has a block structure
    {
    \footnotesize
    \begin{equation*}
    \underbrace{E_1 , E_2, \dots, E_{\sigma_0}}
    _{\text{block 0}}
    ; \,
    \underbrace{E_1(1) , E_2(1), \dots, E_{\sigma_1}(1)}
    _{\text{block 1}}
    ; \,
    \dots; \,
    \underbrace{E_1(m-1) , E_2(m-1), \dots, E_{\sigma_{m-1}}(m-1)}
    _{\text{block $m-1$}}
    \end{equation*}
    }
    where $\sigma = (\sigma_0 \ge \sigma_1 \ge \dots \ge \sigma_{m-1} \ge 0$) is a
    non-increasing sequence of non-negative integers called the \emph{support partition}
    of the collection.
    In the above notation we use semicolons to separate the blocks.
    The block $(E_1 , E_2, \dots, E_{\sigma_0})$ is called the \emph{starting block}.
    We use notation $(E_\bullet,\sigma)$ for a Lefschetz collection with support partition $\sigma$.

    \smallskip

    \item[(ii)] If $\sigma_0 = \sigma_1 = \dots = \sigma_{m-1}$, then the
    Lefschetz collection is called \emph{rectangular}.
    Otherwise, its \emph{rectangular part} is defined to be the subcollection
    \begin{multline}
      E_1, E_2, \dots, E_{\sigma_{m-1}}; E_1(1), E_2(1), \dots, E_{\sigma_{m-1}}(1); \dots \\
      \dots ; E_1(m-1), E_2(m-1), \dots, E_{\sigma_{m-1}}(m-1).
    \end{multline}

    \smallskip

    \item[(iii)] The \emph{residual category} of a Lefschetz collection is defined as the
    orthogonal of its rectangular part
    \begin{equation*}
    \cR =
    \Big\langle
    E_1, E_2, \dots, E_{\sigma_{m-1}};\dots ; E_1(m-1), E_2(m-1), \dots, E_{\sigma_{m-1}}(m-1) \Big\rangle^\perp.
    \end{equation*}
  \end{enumerate}
\end{definition}

\begin{remark}
  \
  \begin{enumerate}
    \item The residual category $\cR$ of a Lefschetz collection vanishes if and only if
    the Lefschetz collection is full and rectangular.

    \item The definitions imply that we have a \emph{semiorthogonal decomposition}
    \begin{equation*}
      \Db(X) = \big\langle \cR ; E_1, E_2, \dots, E_{\sigma_{m-1}};\dots ; E_1(m-1), E_2(m-1), \dots, E_{\sigma_{m-1}}(m-1) \big\rangle.
    \end{equation*}

    \item An important feature of the residual category is the existence of a natural autoequivalence
    $\tau_{\cR} \colon \cR \to \cR$ called the \emph{induced polarization}
    such that $\tau_{\cR}^m \cong \bS_\cR^{-1}[\dim X]$, where~$\bS_\cR$
    is the Serre functor of~$\cR$. It is useful to think of $\tau_{\cR}$
    as the analogue of the twist by $\cO(1)$ in $\Db(X)$.
    For more details on $\tau_{\cR}$ we refer to \cite{KuSm20, KuSm21}.
  \end{enumerate}
\end{remark}

\begin{example}
  \
  \begin{enumerate}
    \item Collection \eqref{eq:collection-beilinson} is a Lefschetz collection with
    \begin{equation*}
      E_\bullet = (\cO) \quad \text{and} \quad \sigma = (1^n) \coloneqq
      (\underbrace{1, \dots, 1}_{\text{$n$ times}}).
    \end{equation*}

    \smallskip

    \item Collection \eqref{eq:collection-kapranov} is a Lefschetz collection with
    \begin{equation*}
      E_\bullet = (\cO, \cU^\vee, S^2\cU^\vee) \quad \text{and} \quad \sigma = (3,2,1).
    \end{equation*}

    \smallskip

    \item On $\G(2,4)$ there exist a Lefschetz collection with a smaller starting
    block than \eqref{eq:collection-kapranov}. Namely, we have
    \begin{equation}\label{eq:collection-minimal}
      \Db(\G(2,4)) = \langle \cO, \cU^\vee; \cO(1), \cU^\vee(1); \cO(2); \cO(3) \rangle,
    \end{equation}
    For this collection we have
    \begin{equation*}
      E_\bullet = (\cO, \cU^\vee) \quad \text{and} \quad \sigma = (2,2,1,1).
    \end{equation*}
  \end{enumerate}
\end{example}

\subsection{Lefschetz collections and quantum cohomology}

Let $X$ now be a Fano variety over~$\Bbbk$.
Roughly speaking the main conjectures of \cite{KuSm20, KuSm21}
say that there is a deep relation between the small quantum cohomology
of $X$ and the structure of Lefschetz collections and their residual categories.
Let us now be more precise.

Let $X$ be a Fano variety with vanishing odd cohomology.
We denote by $\QHcan(X)$ its small quantum cohomology specialized at the canonical class.
If the Picard rank of $X$ is one, and consequently there is only one deformation
parameter $q$ in the small quantum cohomology, $\QHcan(X)$ is the small quantum
cohomology at $q=1$; this is the case for almost all examples appearing below.
Thus, $\QHcan(X)$ is a finite dimensional commutative $\bC$-algebra, whose underlying
vector space is canonically isomorphic to $H^*(X, \bC)$.

Now we define the \textsf{quantum spectrum} of $X$ as
\begin{equation*}
  \QS_{X} \coloneqq \Spec(\QHcan(X)),
\end{equation*}
which is a finite scheme endowed with an action of the group $\mu_m$, where
$m$ is the Fano index of $X$. The anticanonical class $-K_X$ defines a morphism
\begin{equation*}
  \kappa \colon \QS_{X} \to \bA^1,
\end{equation*}
which is equivariant with respect to the standard action of $\mu_m$ on $\bA^1$.
Finally, we define
\begin{equation*}
  \QSx_X \coloneqq \kappa^{-1}(\bA^1 \setminus \{ 0 \}) \quad \text{and} \quad
  \QSo_X \coloneqq \QS_X \setminus \QSx_X.
\end{equation*}
The action of $\mu_m$ on $\QSx_X$ is free, as it is free on $\bA^1 \setminus \{ 0 \}$.
We refer to \cite[Introduction]{KuSm21} for more details on the setup.

\begin{conjecture}[{\cite[Conjecture 1.3]{KuSm21}}]
  \label{conjecture:main}
  Let $X$ be a Fano variety of index~$m$ over an algebraically closed field $\Bbbk$
  of characteristic zero and assume that the big quantum cohomology~$\BQH(X)$ is
  generically semisimple.

  \begin{enumerate}

  \smallskip

  \item There is an $\Aut(X)$-invariant exceptional collection $E_1,\dots,E_k$ in $\Db(X)$, where~$k$
  is the length of~$\QSx_X$ divided by~$m$.
  This collection extends to a rectangular Lefschetz collection
  \begin{multline}\label{eq:main-conjecture-rectangular-part}
  E_1, E_2, \dots, E_k; E_1(1), E_2(1), \dots, E_k(1); \dots \\
  \dots ; E_1(m-1), E_2(m-1), \dots, E_k(m-1).
  \end{multline}
  in $\Db(X)$.

  \medskip

  \item The residual category $\cR$ of \eqref{eq:main-conjecture-rectangular-part}
  has a completely orthogonal $\Aut(X)$-invariant
  decomposition
  \begin{equation*}
  \cR = \bigoplus_{\xi \in \QSo_X} \cR_\xi
  \end{equation*}
  with components indexed by closed points $\xi \in \QSo_X$.
  Moreover, the component~$\cR_\xi$ of~$\cR$ is generated by an exceptional collection
  of length equal to the length of the localization~$(\QSo_X)_\xi$ at~$\xi$.

  \medskip

  \item The induced polarization $\tau_\cR$ permutes the components $\cR_\xi$.
  More precisely, for each point $\xi \in \QSo_X$ it induces an equivalence
  \begin{equation*}
  \tau_\cR \colon \cR_\xi \xrightarrow{\ \sim\ } \cR_{g(\xi)},
  \end{equation*}
  where $g$ is a generator of $\mu_m$.
  \end{enumerate}
\end{conjecture}
Thus, intuitively the points of $\QSx_X$ correspond to the rectangular part
\eqref{eq:main-conjecture-rectangular-part} and the twist by $\cO(1)$ corresponds
to the action of $\mu_m$ on $\QSx_X$;
points of $\QSo_X$ give rise to an exceptional collection in $\cR$.

\medskip

Below we discuss two particular instances of the above conjecture. In the first case
we assume that the small quantum cohomology, or rather $\QHcan(X)$, is semisimple,
and in the second case we consider a particular class of homogeneous varieties
called coadjoint varieties, whose $\QHcan(X)$ is almost never semisimple.

\subsection{Cases with semisimple $\QHcan(X)$}

The simplest example, where Conjecture~\ref{conjecture:main} holds, is provided
by $\bP^n$. Indeed, it is well-known that we have
\begin{equation*}
  \QHcan(\bP^n) = \bC[h]/(h^{n+1}-1),
\end{equation*}
and, therefore, the quantum spectrum $\QS_{\bP^n}$ is a reduced subscheme of $\bA^1$ supported
at the points~$\zeta^i$ with $i \in [0,n]$, where $\zeta$ is a primitive $(n+1)$-st
root of unity.
The action of $\mu_{n+1}$ on $\QS_{\bP^n} = \QSx_{\bP^n}$ is the usual action of $\mu_{n+1}$
on $(n+1)$-st roots of unity. Therefore, this action has only one orbit and,
according to Conjecture \ref{conjecture:main}, in $\Db(\bP^n)$ we should
expect to have a Lefschetz collection, whose starting block $E_\bullet$ consists
of one object and whose support partition is of the form $\sigma = (1^{n+1})$.
Since $\QSo_{\bP^n} = \emptyset$, the residual category vanishes.
The Beilinson collection \eqref{eq:collection-beilinson} satisfies all these requirements.

\smallskip

More generally, if $\QHcan(X)$ is semisimple, then Conjecture \ref{conjecture:main}
gives a full description of the residual category. Indeed, since each component
$\cR_\xi$ is generated by one exceptional object Conjecture \ref{conjecture:main}(ii)
says that the residual category $\cR$ is generated by a completely orthogonal
exceptional collection (cf. \cite[Conjecture 1.12]{KuSm20}).

\smallskip

There is a number of cases with semisimple $\QHcan(X)$, where our conjecture is
known to hold:
\begin{enumerate}
  \item for $\G(k,n)$ with for $k = p$ a prime number (see \cite{KuSm20});

  \item for quadrics $\Q_n$ this follows from Kapranov's work \cite{Ka} (see \cite[Example 1.6]{KuSm20});

  \item for $\OG(2,2n+1)$ this follows from Kuznetsov's work \cite{Ku08a} (see \cite[Example 1.9]{KuSm20});

  \item for $\IG(3,8)$ and $\IG(3,10)$ this holds by \cite{Gu18, No};

  \item for $\IG(4,8)$ and $\IG(5,10)$ this should follow from \cite{Fo19,PS}.

  \item for the Cayley plane $\rE_6/\rP_1$ this holds by \cite{FM, Manivel, BKS};

  \item for the Cayley Grassmannian this holds by \cite{BM, Guseva};

  \item for the $\rG_2$-Grassmannian $\rG_2/\rP_2$ this holds by \cite{K06};

  \item for some horospherical varieties of Picard rank one by \cite{Ku08a, Pech13, GPPS, Fo22};
\end{enumerate}
In the above list $\IG(k,2n)$ is the variety parametrizing $k$-dimensional isotropic
subspaces in a $2n$-dimensional vector space with a symplectic form; this is a
homogeneous space for the symplectic group $\Sp_{2n}$. The variety $\OG(2,2n+1)$
parametrizes $2$-dimensional isotropic subspaces in a $(2n+1)$-dimensional vector
space with a symmetric non-degenerate form; this is a homogeneous space for the
special orthogonal group $\SO_{2n+1}$. We refrain from recalling the definitions
of all the other varieties in the above list.

\subsection{Cases with non-semisimple $\QHcan(X)$}

If the algebra $\QHcan(X)$ is not semisimple, then Conjecture \ref{conjecture:main}
does not give a full description of the orthogonal components~$\cR_\xi$ of the residual category.
However, under mirror symmetry, the locus $\QSo_X$ corresponds to the critical points
of the mirror Landau--Ginzburg model $f$ of the Fano variety~$X$ in the fiber over
zero $f^{-1}(0)$.
We expect that for each $\xi \in \QSo_X$ the component $\cR_\xi$ is equivalent
to the Fukaya--Seidel category of the corresponding critical point in $f^{-1}(0)$.
Below we illustrate this phenomenon.

\medskip
\noindent
\textbf{Symplectic isotropic Grassmannians $\IG(2,2n)$.}
The variety $\IG(2,2n)$ parametrizes $2$-dimensional isotropic subspaces in a
$2n$-dimensional vector space with a symplectic form. It is embedded into the
usual Grassmannian $\IG(2,2n) \subset \G(2,2n)$ and is, in fact, a hyperplane
section of $\G(2,2n)$.

The Fano index of $\IG(2,2n)$ is equal to $2n-1$ and
$\dim_{\bC}(H^{*}(\IG(2,2n), \bC)) = 2n(n-1)$.
We know by \cite[Proposition 4.3]{CMKMPS}
that $\QSx_{\IG(2,2n)}$ is a disjoint union of $n-1$ orbits of $\mu_{2n-1}$,
each of which consists of $2n-1$ reduced points, and $\QSo_{\IG(2,2n)}$ consists of one
non-reduced point $\xi_0$ with local algebra $\bC[t]/t^{n-1}$.

Since $\bC[t]/t^{n-1}$ is isomorphic to the Jacobi algebra of an isolated
hypersurface singularity of type $\rA_{n-1}$, the above description of $\QSo_{\IG(2,2n)}$
suggests that in Conjecture \ref{conjecture:main} we should expect the residual
category $\cR = \cR_{\xi_0}$ to be equivalent to the Fukaya--Seidel category of
an isolated hypersurface singularity of type $\rA_{n-1}$, which by \cite{Seidel}
is equivalent to the derived category of representations of a quiver of type $\rA_{n-1}$.

A Lefschetz collection on $\IG(2,2n)$ satisfying these conjectures was constructed in \cite{Ku08a}.
Indeed, let us define
\begin{equation*}
  \begin{aligned}
    & E_i = S^{i-1} \cU^\vee \quad \text{for} \quad 1 \leq i \leq n, \\
    & \sigma = (n^{n-1}, (n-1)^n).
  \end{aligned}
\end{equation*}
Then $(E_\bullet, \sigma)$ is a full Lefschetz collection by \cite[Theorem 5.1]{Ku08a}.
Moreover, by \cite[Theorem 9.6]{CMKMPS}, its residual category is equivalent
to the derived category of representations of $\rA_{n-1}$ quiver.

\medskip
\noindent
\textbf{Coadjoint varieties.}
The variety $\IG(2,2n)$ considered above fits naturally into a series of examples.
Namely, $\IG(2,2n)$ is the coadjoint variety in Dynkin type $\rC_n$.
In general, the \emph{coadjoint variety} of a simple algebraic group $\rG$ is
the highest weight vector orbit in the projectivization of the irreducible
$\rG$-representation, whose highest weight is the highest short root.
Therefore, a coadjoint variety is uniquely determined by the Dynkin type of $\rG$,
which can be $\rA_n, \rB_n, \rC_n, \rD_n, \rE_6, \rE_7, \rE_8, \rF_4, \rG_2$.
Coadjoint varieties are of Picard rank one except for type $\rA_n$, where the Picard rank is two.

The quantum spectrum of a coadjoint variety is described in \cite{PeSm}. To state
the result we denote by $\rT(\rG)$ the Dynkin diagram of $\rG$ and by $\rTs(\rG)$
the subdiagram of $\rT(\rG)$ consisting of vertices corresponding to short roots.
Thus, we have the table
\begin{equation*}
\begin{array}{|c|c|c|c|c|c|c|c|}
\hline
\rT & \rA_n & \rB_n & \rC_n & \rD_n & \rE_n & \rF_4 & \rG_2 \\
\hline
\rTs & \rA_n & \rA_1 & \rA_{n-1} & \rD_n & \rE_n & \rA_2 & \rA_1 \\
\hline
\end{array}
\end{equation*}

Now we are ready to describe the quantum spectrum of coadjoint varieties.
\begin{lemma}[{\cite{PeSm}}]
  Let $X$ be the coadjoint variety of a simple algebraic group~$\rG$. Then
  \begin{enumerate}
    \item $\QSx(X)$ consists of reduced points;

    \item if $\rT(\rG) = \rA_n$ and $n$ is even, then $\QSo(X) = \emptyset$;

    \item otherwise, $\QSo(X)$ has a unique point and the local algebra at this
    point is isomorphic to the Jacobi ring of an isolated hypersurface singularity
    of Dynkin type $\rTs(\rG)$.
  \end{enumerate}
\end{lemma}

Similarly to the case of $\IG(2,2n)$ the above description suggests the following.

\begin{conjecture}[{\cite[Conjecture 1.8]{KuSm21}}]
\label{conjecture:coadjoint}
Let $X$ be the coadjoint variety of a simple algebraic group~$\rG$ over an algebraically closed field of characteristic zero.
Then $\Db(X)$ has an~$\Aut(X)$-invariant rectangular Lefschetz exceptional collection with residual category~$\cR$ and
\begin{enumerate}
\item
if $\rT(\rG) = \rA_n$ and $n$ is even, then $\cR = 0$;
\item
otherwise, $\cR$ is equivalent to the derived category of representations of a quiver of Dynkin type~$\rTs(\rG)$.
\end{enumerate}
\end{conjecture}

This conjecture is by now known in all Dynkin types except for $\rE_6, \rE_7, \rE_8$.
Indeed, this is
\cite[Theorem~9.6]{CMKMPS} for type~$\rC_n$,
\cite[Example~1.6]{KuSm20} for types $\rB_n$ and~$\rG_2$,
\cite[Theorem 1.9]{KuSm21} for types $\rA_n$ and $\rD_n$,
\cite[Theorem~1.4]{BKS} for type~$\rF_4$.

\begin{remark}
  Above we have discussed only coadjoint varieties, but one can also define adjoint varieties.
  The \emph{adjoint variety} of a simple algebraic group $\rG$ is
  the highest weight vector orbit in the projectivization of the irreducible
  $\rG$-representation, whose highest weight is the highest root. If $\rT(\rG)$
  is simply laced, then all roots have the same length and, therefore, adjoint
  and coadjoint varieties coincide.

  Let $X$ be the adjoint variety of a simple algebraic group $\rG$ whose Dynkin type
  is not simply laced, i.e. $\rB_n, \rC_n, \rF_4, \rG_2$. Then by \cite[Theorem 9.1]{PeSm}
  we know that $\QSo(X) = \emptyset$ and by Conjecture \ref{conjecture:main} we should
  expect a rectangular Lefschetz collection in $\Db(X)$.
  This is by now known in all cases:
  \cite[Theorem~7.1]{Ku08a} for type~$\rB_n$,
  \cite[Example~1.4]{KuSm20} for type~$\rC_n$,
  \cite[Theorem~1.1]{Sm21} for type~$\rF_4$,
  and \cite[\S6.4]{K06} for type~$\rG_2$.
\end{remark}

\medskip
\noindent
\textbf{Acknowledgements.}
I am very grateful to my coauthor Alexander Kuznetsov for the joint projects
\cites{KuSm20, KuSm21} on which this talk is based and for the comments on the
preliminary version of this note.

Finally, I thank the organisers of the Nottingham Online Algebraic Geometry Seminar,
and in particular Al Kasprzyk, for the opportunity to present this work at the seminar.

\vspace{10pt}

\bibliographystyle{plain}
\bibliography{refs}

\begin{bibdiv}
\begin{biblist}

\bib{Be}{article}{
      author={Beilinson, A.~A.},
       title={Coherent sheaves on {${\bf P}\sp{n}$}\ and problems in linear
  algebra},
        date={1978},
        ISSN={0374-1990},
     journal={Funktsional. Anal. i Prilozhen.},
      volume={12},
      number={3},
       pages={68\ndash 69},
}

\bib{BKS}{article}{
      author={Belmans, P.},
      author={Kuznetsov, A.},
      author={Smirnov, M.},
       title={Derived categories of the {C}ayley plane and the coadjoint
  {G}rassmannian of type~$\mathrm{F}$},
        date={2021},
     journal={Transformation Groups},
      eprint={https://doi.org/10.1007/s00031-021-09657-w},
         url={https://doi.org/10.1007/s00031-021-09657-w},
        note={See also
  \href{https://arxiv.org/abs/2005.01989}{arXiv:2005.01989}},
}

\bib{BM}{article}{
      author={Benedetti, V.},
      author={Manivel, L.},
       title={The small quantum cohomology of the {C}ayley {G}rassmannian},
        date={2020},
     journal={Internat. J. Math.},
      volume={31},
      number={3},
       pages={2050019, 23},
         url={https://doi.org/10.1142/S0129167X20500196},
}

\bib{CMKMPS}{article}{
      author={Cruz~Morales, J.~A.},
      author={Kuznetsov, A.},
      author={Mellit, A.},
      author={Perrin, N.},
      author={Smirnov, M.},
       title={On quantum cohomology of {Grassmannian}s of isotropic lines,
  unfoldings of {$A_n$}-singularities, and {Lefschetz} exceptional
  collections},
        date={2019},
     journal={Annales de l'Institut Fourier},
      volume={69},
      number={3},
       pages={955\ndash 991},
        note={See also
  \href{https://arxiv.org/abs/1705.01819}{arXiv:1705.01819}},
}

\bib{FM}{article}{
      author={Faenzi, D.},
      author={Manivel, L.},
       title={On the derived category of the {C}ayley plane {II}},
        date={2015},
     journal={Proc. Amer. Math. Soc.},
      volume={143},
      number={3},
       pages={1057\ndash 1074},
        note={See also
  \href{https://arxiv.org/abs/1201.6327}{arXiv:1201.6327}},
}

\bib{Fo19}{article}{
      author={Fonarev, A.},
       title={Full exceptional collections on {L}agrangian {G}rassmannians},
        date={2020},
     journal={Int. Math. Res. Not.},
        note={Available online ahead of print as
  \href{https://doi.org/10.1093/imrn/rnaa098}{doi:10.1093/imrn/rnaa098}, see
  also \href{https://arxiv.org/abs/1911.08968}{arXiv:1911.08968}},
}

\bib{Fo22}{article}{
      author={Fonarev, A.},
       title={On the bounded derived category of {IG}r(3, 7)},
        date={2022},
     journal={Transform. Groups},
      volume={27},
      number={1},
       pages={89\ndash 112},
         url={https://doi.org/10.1007/s00031-020-09614-z},
}

\bib{GPPS}{article}{
      author={Gonzales, R.},
      author={Pech, C.},
      author={Perrin, N.},
      author={Samokhin, A.},
       title={Geometry of horospherical varieties of {P}icard rank one},
        note={Available at
  \href{https://arxiv.org/abs/1803.05063}{arXiv:1803.05063}},
}

\bib{Guseva}{article}{
      author={Guseva, L.},
       title={{O}n the derived category of the {C}ayley {G}rassmannian},
        note={to appear},
}

\bib{Gu18}{article}{
      author={Guseva, L.},
       title={On the derived category of {$\mathrm{IGr}(3;8)$}},
        date={2020},
        ISSN={0368-8666},
     journal={Mat. Sb.},
      volume={211},
      number={7},
       pages={24\ndash 59},
        note={See also
  \href{https://arxiv.org/abs/1810.07777}{arXiv:1810.07777}},
}

\bib{Ka}{article}{
      author={Kapranov, M.~M.},
       title={On the derived categories of coherent sheaves on some homogeneous
  spaces},
        date={1988},
     journal={Inventiones Mathematicae},
      volume={92},
      number={3},
       pages={479\ndash 508},
}

\bib{K06}{article}{
      author={Kuznetsov, A.},
       title={Hyperplane sections and derived categories},
        date={2006},
     journal={Izv. Ross. Akad. Nauk Ser. Mat.},
      volume={70},
      number={3},
       pages={23\ndash 128},
        note={See also
  \href{https://arxiv.org/abs/math/0503700}{arXiv:math/0503700}},
}

\bib{K07}{article}{
      author={Kuznetsov, A.},
       title={Homological projective duality},
        date={2007},
        ISSN={0073-8301},
     journal={Publ. Math. Inst. Hautes \'{E}tudes Sci.},
      number={105},
       pages={157\ndash 220},
         url={https://doi.org/10.1007/s10240-007-0006-8},
        note={See also
  \href{https://arxiv.org/abs/math/0507292}{arXiv:math/0507282}},
}

\bib{Ku08a}{article}{
      author={Kuznetsov, A.},
       title={Exceptional collections for {G}rassmannians of isotropic lines},
        date={2008},
     journal={Proceedings of the London Mathematical Society. Third Series},
      volume={97},
      number={1},
       pages={155\ndash 182},
        note={See also
  \href{https://arxiv.org/abs/math/0512013}{arXiv:math/0512013}},
}

\bib{KuSm20}{article}{
      author={Kuznetsov, A.},
      author={Smirnov, M.},
       title={On residual categories for {G}rassmannians},
        date={2020},
     journal={Proceedings of the London Mathematical Society},
      volume={120},
      number={5},
       pages={617\ndash 641},
}

\bib{KuSm21}{article}{
      author={Kuznetsov, A.},
      author={Smirnov, M.},
       title={Residual categories for (co)adjoint {G}rassmannians in classical
  types},
        date={2021},
     journal={Compositio Mathematica},
      volume={157},
      number={6},
       pages={1172–1206},
        note={See also
  \href{https://arxiv.org/abs/2001.04148}{arXiv:2001.04148}},
}

\bib{Manivel}{article}{
      author={Manivel, L.},
       title={On the derived category of the {C}ayley plane},
        date={2011},
     journal={J. Algebra},
      volume={330},
       pages={177\ndash 187},
        note={See also
  \href{https://arxiv.org/abs/0907.2784}{arXiv:0907.2784}},
}

\bib{No}{article}{
      author={Novikov, A.},
       title={{L}efcshetz {E}xceptional {C}ollections on {I}sotropic
  {G}rassmannians},
        date={2020},
        note={Master thesis at the Higher School of Economics, Moscow.
  Available at {{\tt https://www.hse.ru/en/edu/vkr/369850047}}},
}

\bib{Pech13}{article}{
      author={Pech, C.},
       title={Quantum cohomology of the odd symplectic {G}rassmannian of
  lines},
        date={2013},
     journal={J. Algebra},
      volume={375},
       pages={188\ndash 215},
         url={https://doi.org/10.1016/j.jalgebra.2012.11.010},
}

\bib{PeSm}{article}{
      author={Perrin, N.},
      author={Smirnov, M.},
       title={On the big quantum cohomology of coadjoint varieties},
        note={Available at
  \href{https://arxiv.org/abs/2112.12436}{arXiv:2112.12436}},
}

\bib{PS}{article}{
      author={Polishchuk, A.},
      author={Samokhin, A.},
       title={Full exceptional collections on the {L}agrangian {G}rassmannians
  {$LG(4,8)$} and {$LG(5,10)$}},
        date={2011},
     journal={Journal of Geometry and Physics},
      volume={61},
      number={10},
       pages={1996\ndash 2014},
}

\bib{Seidel}{incollection}{
      author={Seidel, P.},
       title={More about vanishing cycles and mutation},
        date={2001},
   booktitle={Symplectic geometry and mirror symmetry ({S}eoul, 2000)},
   publisher={World Sci. Publ., River Edge, NJ},
       pages={429\ndash 465},
}

\bib{Sm21}{article}{
      author={Smirnov, M.},
       title={On the derived category of the adjoint {G}rassmannian of type
  {F}},
        note={Available at
  \href{https://arxiv.org/abs/2107.07814}{arXiv:2107.07814}},
}

\end{biblist}
\end{bibdiv}

\end{document}